\documentclass{amsart}
\usepackage{amsmath,amsthm,amssymb,amsfonts}

\theoremstyle{plain}
\newtheorem{theo+}           {Theorem}      [section]
\newtheorem{prop+}  [theo+]  {Proposition}
\newtheorem{lemm+}  [theo+]  {Lemma}

\newenvironment{theorem}{\begin{theo+}}{\end{theo+}}
\newenvironment{proposition}{\begin{prop+}}{\end{prop+}}
\newenvironment{lemma}{\begin{lemm+}}{\end{lemm+}}

\newcommand{\Z}{{\mathbb Z}}
\newcommand{\Zp}{{\mathbb Z}_{\geq 0}}
\newcommand{\Zn}{{\mathbb Z}_{< 0}}
\newcommand{\la}{\lambda}
\newcommand{\La}{\Lambda}
\newcommand{\om}{\omega}
\newcommand{\Om}{\Omega}

\begin{document}
\baselineskip 18pt
\larger[2]
\title[A multivariable elliptic summation formula]
{A proof of a multivariable elliptic
summation\\ formula conjectured by Warnaar}
\author{Hjalmar Rosengren}
\address{Department of Mathematics\\ Chalmers University of Technology and 
G\"oteborg University\\SE-412~96 G\"oteborg, Sweden}
\email{hjalmar@math.chalmers.se}

\begin{abstract}
We prove a multivariable elliptic analogue of Jackson's ${}_8W_7$
summation formula, which was recently conjectured by S.\ O.\ Warnaar. 
\end{abstract}
\maketitle   

\section{Introduction}
Elliptic hypergeometric series form a natural 
 generalization of hypergeometric and basic hypergeometric (or $q$-) series.
It is surprising that they were introduced only very  recently,
by Frenkel and Turaev \cite{ft}, who expressed the
$6j$-symbols corresponding to certain elliptic solutions of the 
Yang--Baxter equation, cf.\ \cite{dj},  in terms of the
 ${}_{10}\om_9$-sums defined below. 
It is expected that elliptic hypergeometric series
play a fundamental role in the representation theory of elliptic
quantum groups, though so far there has been little
work in this direction. 

Recall that a series  $\sum_{n} a_n$ is called hypergeometric if
 $f(n)=a_{n+1}/a_n$ is a rational function of $n$ and basic hypergeometric
if $f$ is a rational function of $q^n$ for some $q$.
 This can be compared with Weierstrass'
 theorem, stating that a meromorphic function of $z$ which satisfies an
algebraic addition theorem is either a rational function, a rational
function of $q^z$, or, in the most general case,
 an elliptic function. This suggests that an elliptic hypergeometric
series should be a series $\sum_{n} a_n$
with $a_{n+1}/a_n$  an elliptic function of $n$.
Actually, the series introduced by Frenkel and Turaev only fit
this description if one interprets the term ``elliptic'' somewhat
loosely. Nevertheless, their properties stem from addition
theorems for elliptic functions
  (it is worth noting that the  Yang--Baxter equation is an algebraic
addition theorem for matrix-valued functions).

Let us write $[x]$ for the ``elliptic number''
(a Jacobi theta function, normalized so that $[1]=1$)
$$[x]=\frac{q^{-\frac x2}\prod_{j=0}^\infty(1-q^xp^j)(1-q^{-x}p^{j+1})}
{q^{-\frac 12}\prod_{j=0}^\infty(1-qp^j)(1-q^{-1}p^{j+1})},$$
 where $p$ and $q$ are fixed parameters with $|p|<1$.
When $p=0$, $q=e^{2ih}$, we have the trigonometric number
$$[x]=\frac{q^{\frac x2}-q^{-\frac x2}}{q^{\frac 12}-q^{-\frac 12}}
=\frac{\sin(hx)}{\sin(h)},$$
which tends to the rational number $[x]=x$ as $q$ tends to $1$.
Returning to the general case, we write
$$[x]_n=[x][x+1]\dotsm [x+n-1]
$$
for the elliptic Pochhammer symbols.
The elliptic, or modular, hypergeometric series
 occurring in \cite{ft} are finite sums of the form
\begin{equation*}\begin{split}&\quad{}_{r+1}\om_r(a;-n,b_1,\dots,b_{r-3})\\
&= \sum_{k=0}^n \frac{[a+2k]}{[a]}
\frac{[a]_k[-n]_k[b_1]_k\dotsm[b_{r-3}]_k}
{[1]_k[1+a+n]_k[1+a-b_1]_k\dotsm[1+a-b_{r-3}]_k},\end{split}\end{equation*}
where
\begin{equation}\label{bal}(r-3)(a+1)=2\left(1-n+\textstyle\sum_i\displaystyle
 b_i\right).
\end{equation}
When $p=0$, this is a terminating
very-well-poised balanced basic hypergeometric series
\cite{gr}, which tends to the corresponding
  hypergeometric series as $q$ tends to $1$.
As was pointed out in \cite{ft}, the series ${}_{r+1}\om_r$
has remarkable invariance 
properties under the standard action of $\mathrm{SL}(2,\mathbb Z)$
on $p$ and $q$.

Most (or possibly all) known identities involving terminating $q$-series
 may be proved by induction, using the trigonometric
addition formula
\begin{equation}\label{tat}
[x+z]_{p=0}[x-z]_{p=0}=[x+y]_{p=0}[x-y]_{p=0}+[y+z]_{p=0}[y-z]_{p=0}.
\end{equation}
However, only a tiny subset of these identities may be obtained from the
elliptic addition formula
$$[x+z][x-z][y+w][y-w]=[x+y][x-y][z+w][z-w]+[x+w][x-w][y+z][y-z]$$
satisfied by the  elliptic numbers. 
At least as a rule of thumb,  these are 
the identities involving series which are
both well-poised and balanced, and  thus only these admit elliptic
analogues. In particular, Frenkel and Turaev obtained the elliptic
  Jackson--Dougall summation formula
\begin{equation}\label{ej}{}_8\om_7(a;-n,b,c,d,e)
=\frac{[a+1]_n[a+1-b-c]_n[a+1-b-d]_n[a+1-c-d]_n}
{[a+1-b]_n[a+1-c]_n[a+1-d]_n[a+1-b-c-d]_n}\end{equation}
and (more generally) the elliptic Bailey transformation
formula  
\begin{equation*}\begin{split}{}_{10}\om_9(a;-n,b,c,d,e,f,g)&=
\frac{[a+1]_n[a+1-e-f]_n[\la+1-e]_n[\la+1-f]_n}
{[a+1-e]_n[a+1-f]_n[\la+1-e-f]_n[\la+1]_n}\\
&\quad\times\,{}_{10}\om_9(\la;-n,\la+b-a,\la+c-a,\la+d-a,e,f,g),
\end{split}\end{equation*}
where $\la=2a+1-b-c-d$;
note that the balanced condition \eqref{bal} is assumed.

If one wants to further develop the theory of elliptic hypergeometric series,
there are two natural directions: quadratic (or higher) transformation
formulas and multivariable series. 
In \cite{w}, Warnaar initiated the investigation of both topics.
We will be concerned with the multivariable theory. As Warnaar pointed out,
progress in this direction requires essentially new ideas,
since the known proofs in the trigonometric and rational case usually
depend  on ``lower level'' identities, corresponding
to the degenerate addition theorem \eqref{tat}. 

The purpose of this paper is to prove an identity conjectured by  
Warnaar in \cite{w},
cf.\ Theorem \ref{t}, which is a
generalization of \eqref{ej} connected with the root system $C_n$.
Our main tool will be a different generalization of \eqref{ej}, 
obtained by Warnaar \cite{w} from a determinant evaluation. 

We mention that one degenerate case of Theorem \ref{t}
is the terminating case of a multivariable ${}_6\psi_6$ sum
due to van Diejen \cite{vd}. 
It generalizes various Macdonald--Morris-type identities for root systems,
cf.\ \cite{vd} for a detailed discussion.
Moreover, van Diejen's sum
gives the norm evaluation for the multivariable $q$-Racah polynomials  
studied by van Diejen and Stokman \cite{dst}. 

When \cite{w} was published, Theorem \ref{t} was new even in the
trigonometric case ($p=0$). This case of the conjecture was settled 
by van Diejen and Spiridonov \cite{ds}, who deduced it from a certain
multiple integral due to Gustafson \cite{g}, 
which reduces to  the
Nassrallah--Rahman integral \cite{nr} in the one-variable case. The multiple
$q$-series in question appears as a sum of residues of the integrand.
Moreover, it was demonstrated that both sides of the 
 equality in Theorem \ref{t} are invariant under the action of
$\mathrm{SL}(2,\mathbb Z)$. Using the theory of modular forms, 
it was then proved that for $q=e^{2 i h}$, the two sides are equal 
at least up to order $h^{10}$ around $h=0$; a strong indication that
Warnaar's conjecture is true. Finally,
van Diejen and Spiridonov conjectured an elliptic generalization of 
Gustafson's
integral, involving the elliptic gamma function introduced by
Ruijsenaars \cite{r}. A  proof of this identity would yield
 another proof of Theorem \ref{t},  completely
different from the one given here. 
The one-variable case of the integral is treated in \cite{sp}.

{\bf Acknowledgement:} I would like to thank Jan Felipe van Diejen and
 Vyacheslav Spiridonov for illuminating correspondence.

\section{Notation and statement of results}
In the rest of the paper
we will use the ``multiplicative'' notation of \cite{w}
rather than the ``additive'' notation of \cite{ft}
used in the introduction.
Since the elliptic modulus $p$ is 
fixed we suppress it from the notation. Thus we  write
$$E(x)=\prod_{j=0}^\infty \left(1-xp^j\right)\left(1-p^{j+1}/x\right),$$
$$E(x_1,\dots,x_m)=E(x_1)\dotsm E(x_m),$$
$$(a;q)_k=\prod_{j=0}^{k-1}E(aq^j),\ \ \ k\in\Zp,$$
$$(a;q)_{-k}=\frac{1}{(aq^{-k};q)_k},\ \ \ k\in\Zn,$$
$$(a_1,\dots,a_m;q)_k=(a_1;q)_k\dotsm(a_m;q)_k, $$
$$(a;q,x)_{\la}=\prod_{j=1}^n(ax^{1-j};q)_{\la_j},\ \ \ \la\in\Z^{n},$$
$$(a_1,\dots,a_m;q,x)_\la=(a_1;q,x)_\la\dotsm(a_m;q,x)_\la.$$
We will use without comment standard identities such as
$(a;q)_n(aq^n;q)_k=(a;q)_{n+k}$. 
 We also mention the easily verified identity
\begin{equation}\label{dp}(aq;q)_n\prod_{1\leq i<j\leq n}E(aq^{i+j})=
(aq;q^2)_n\prod_{1\leq i<j\leq n}E(aq^{i+j-1}).\end{equation}

We can now state the  main result of the paper, conjectured by 
Warnaar \cite{w}.

\begin{theorem}\label{t} In the notation above,
\begin{multline}\label{w}
\sum_{\la}\prod_{i=1}^n \left(\frac{E(ax^{2(1-i)}q^{2\la_i})}{E(ax^{2(1-i)})}
\,q^{\la_i}x^{2(i-1)\la_i}\right)\prod_{1\leq i<j\leq n}
\left(\frac{E(x^{j-i}
q^{\la_i-\la_j})}{E(x^{j-i})}\frac{E(ax^{2-i-j}q^{\la_i+\la_j})}
{E(ax^{2-i-j})}\right.\\
\times\left.\frac{(ax^{3-i-j};q)_{\la_i+\la_j}
(x^{j-i+1};q)_{\la_i-\la_j}}{(aqx^{1-i-j};q)_{\la_i+\la_j}
(qx^{j-i-1};q)_{\la_i-\la_j}}\right)
\frac{(ax^{1-n},b,c,d,e,q^{-N};q,x)_\la}{(qx^{n-1},aq/b,aq/c,aq/d,
aq/e,aq^{N+1};q,x)_\la}\\
=\frac{(aq,aq/bc,aq/bd,aq/cd;q,x)_{N^n}}
{(aq/b,aq/c,aq/d,aq/bcd;q,x)_{N^n}},\hfill 
\end{multline}
 where the sum is over the partitions
$$\la\in\La_{nN}=\left\{\la\in\Z^n;N\geq\la_1\geq\la_2\geq\dots\geq\la_n\geq 0
\right\}$$
and where
$bcdex^{n-1}=a^2q^{N+1}$. 
\end{theorem}

Here $N^n$ denotes the partition with
$\la_i=N$, $i=1,\dots,n$. Our main tool will
be the following identity, again due to Warnaar.

\begin{lemma}\label{l}In the notation above,
\begin{equation}\label{w2}\begin{split}
&\sum_{k_1,\dots,k_n=0}^1\, \prod_{i=1}^n
\frac{(bx_i,cx_i,dx_i,ex_i;q)_{k_i}}
{(aqx_i/b,aqx_i/c,aqx_i/d,aqx_i/e;q)_{k_i}}\,(-1)^{k_i}q^{(i-1)k_i}\\
&\quad\times\prod_{1\leq i< j\leq n} \frac{E(q^{k_i-k_j}x_i/x_j)}
{E(x_i/x_j)}\frac{E(ax_ix_jq^{k_i+k_j})}{E(ax_ix_jq)}\\
&=(aq/bc,aq/bd,aq/cd;q^{-1})_n\prod_{i=1}^n\frac{E(aqx_i^2)}
{E(aq^{2-n}/bcdx_i,aqx_i/b,aqx_i/c,aqx_i/d)},
\end{split}\end{equation}
where $a^2q^{3-n}=bcde$.
\end{lemma}

In fact, Warnaar proved the more general identity
\cite[Theorem 5.1]{w} 
\begin{equation*}\begin{split}
&\sum_{k_1,\dots,k_n=0}^N\, \prod_{i=1}^n\frac{E(ax_i^2q^{2k_i})}{E(ax_i^2)}
\frac{(ax_i^2,bx_i,cx_i,dx_i,ex_i,q^{-N};q)_{k_i}}
{(q,aqx_i/b,aqx_i/c,aqx_i/d,aqx_i/e,aq^{N+1}x_i^2;q)_{k_i}}\,q^{ik_i}\\
&\quad\times\prod_{1\leq i< j\leq n} \frac{E(q^{k_i-k_j}x_i/x_j)}
{E(x_i/x_j)}\frac{E(ax_ix_jq^{k_i+k_j})}{E(ax_ix_jq^N)}\\
&=\prod_{i=1}^n\frac{(aqx_i^2,aq^{2-i}/bc,aq^{2-i}/bd,aq^{2-i}/cd;q)_N}
{(aq^{2-n}/bcdx_i,aqx_i/b,aqx_i/c,aqx_i/d;q)_N},
\end{split}\end{equation*}
where $a^2q^{N+2-n}=bcde$. For $n=1$, this is equivalent to \eqref{ej} 
and for $N=1$ it reduces to \eqref{w2}. The case $p=0$ 
is due to Schlosser \cite{s}.

We will prove Theorem \ref{t} by induction on the ``terminator'' $N$.
However, because of a duality property for the sums in question,
cf.\ Proposition \ref{p}, we can alternatively formulate the proof
as an induction on  the number $n$ of variables. 
In this context we remark that,
for $p=0$, \eqref{w2} is a special case not only
of Schlosser's identity but also of yet another multivariable
Jackson--Dougall formula due to Denis and
Gustafson \cite{dg}
and Milne and Lilly \cite{ml}. The degeneration of the latter
 to the
${}_6\psi_6$-level, together with induction on the number of variables,
 was used by van Diejen 
\cite{vd} to prove the trigonometric ${}_6\psi_6$-version of Theorem \ref{t}.
Nevertheless, our proof  is essentially different from the 
one in \cite{vd}, since we only need a very special case of the 
(as yet unproved)
elliptic Denis--Gustafson--Milne--Lilly identity.

\section{Proof of Theorem \ref{t}}

We will prove Theorem \ref{t} by induction on $N$. The argument allows us
to deduce the case $N=1$, equivalent to \eqref{ej}, from the trivial
case $N=0$, so we obtain in particular a direct proof of the
 one-variable elliptic Jackson--Dougall formula.

Assume that Theorem \ref{t} holds for a fixed 
value  of $N$. Let us  fix parameters with
\begin{equation}\label{bc}bcdex^{n-1}=a^2q^{N+2}.\end{equation}
 We  write the right-hand side of \eqref{w} 
 with $N$ replaced by $N+1$ as  
\begin{equation*}\begin{split}R&=\frac{(aq,aq/bc,aq/bd,aq/cd;q,x)_{(N+1)^n}}
{(aq/b,aq/c,aq/d,aq/bcd;q,x)_{(N+1)^n}}\\
&=\frac{(aq,aq/bc,aq/bd,aq/cd;x^{-1})_n}
{(aq/b,aq/c,aq/d,aq/bcd;x^{-1})_n}
\frac{(aq^2,aq^2/bc,aq^2/bd,aq^2/cd;q,x)_{N^n}}
{(aq^2/b,aq^2/c,aq^2/d,aq^2/bcd;q,x)_{N^n}},\end{split}\end{equation*}
where the second factor is the right-hand side of \eqref{w} with $a$ replaced
by $aq$ and $e$ by $eq$. Using our induction hypothesis, we have 
\begin{multline*}R=
\frac{(aq,aq/bc,aq/bd,aq/cd;x^{-1})_n}
{(aq/b,aq/c,aq/d,aq/bcd;x^{-1})_n}
\sum_{\la}\prod_{i=1}^n \left(\frac{E(ax^{2(1-i)}q^{2\la_i+1})}
{E(ax^{2(1-i)}q)}
\,q^{\la_i}x^{2(i-1)\la_i}\right)\\
\begin{split}&\times\prod_{1\leq i<j\leq n}
\left(\frac{E(x^{j-i}
q^{\la_i-\la_j})}{E(x^{j-i})}\frac{E(ax^{2-i-j}q^{\la_i+\la_j+1})}
{E(ax^{2-i-j}q)}\frac{(aqx^{3-i-j};q)_{\la_i+\la_j}
(x^{j-i+1};q)_{\la_i-\la_j}}{(aq^2x^{1-i-j};q)_{\la_i+\la_j}
(qx^{j-i-1};q)_{\la_i-\la_j}}\right)\\
&\times\frac{(aqx^{1-n},b,c,d,eq,q^{-N};q,x)_\la}
{(qx^{n-1},aq^2/b,aq^2/c,aq^2/d,aq/e,aq^{N+2};q,x)_\la}.
\end{split}\end{multline*}

We now  apply \eqref{w2} with 
$$(a,b,c,d,e,x_i,q)\mapsto(aq/x,b,c,d,eq^{-N},x^{1-i}q^{\la_i},x),$$
 which allows us to write
\begin{equation*}\begin{split}&\quad\quad(aq/bc,aq/bd,aq/cd;x^{-1})_n
\prod_{i=1}^n E(ax^{2(1-i)}q^{2\la_i+1})\\
&=\prod_{i=1}^n E(ax^{i-n}q^{1-\la_i}/bcd,ax^{1-i}q^{\la_i+1}/b,
ax^{1-i}q^{\la_i+1}/c, ax^{1-i}q^{\la_i+1}/d)\\
&\quad\times\sum_{k_1,\dots,k_n=0}^1\, \prod_{i=1}^n
\frac{(bx^{1-i}q^{\la_i},cx^{1-i}q^{\la_i},dx^{1-i}q^{\la_i},
ex^{1-i}q^{\la_i-N};x)_{k_i}(-1)^{k_i}x^{(i-1)k_i}}
{(ax^{1-i}q^{\la_i+1}/b,ax^{1-i}q^{\la_i+1}/c,ax^{1-i}q^{\la_i+1}/d,
ax^{1-i}q^{\la_i+N+1}/e;x)_{k_i}}\\
&\quad\times\prod_{1\leq i< j\leq n}
\frac{E(x^{j-i+k_i-k_j}q^{\la_i-\la_j})}
{E(x^{j-i}q^{\la_i-\la_j})}\frac{E(ax^{1-i-j+k_i+k_j}q^{\la_i+\la_j+1})}
{E(ax^{2-i-j}q^{\la_i+\la_j+1})}.
\end{split}\end{equation*}
Plugging this into the previous identity and then replacing $\la$ by $\la-k$
in the summation yields
\begin{multline*}R=\frac{(aq;x^{-1})_n}
{(aq/b,aq/c,aq/d,aq/bcd;x^{-1})_n}\\
\times
\sum_{\la,k}\prod_{i=1}^n \left(\frac{E(ax^{i-n}q^{k_i-\la_i+1}/bcd,
ax^{1-i}q^{\la_i-k_i+1}/b,ax^{1-i}q^{\la_i-k_i+1}/c,
ax^{1-i}q^{\la_i-k_i+1}/d)}{E(aqx^{2(1-i)})}\right.\hfill\\
\times\left.\frac{(bx^{1-i}q^{\la_i-k_i},cx^{1-i}q^{\la_i-k_i},
dx^{1-i}q^{\la_i-k_i},ex^{1-i}q^{\la_i-k_i-N};x)_{k_i}
(-1)^{k_i}q^{\la_i-k_i}x^{(i-1)(2\la_i-k_i)}}
{(ax^{1-i}q^{\la_i-k_i+1}/b,ax^{1-i}q^{\la_i-k_i+1}/c,
ax^{1-i}q^{\la_i-k_i+1}/d,ax^{1-i}q^{\la_i-k_i+N+1}/e;x)_{k_i}}
\right)\hfill\\
\times\prod_{1\leq i<j\leq n}\left\{\frac{E(x^{j-i+k_i-k_j}
q^{\la_i-\la_j-k_i+k_j})}{E(x^{j-i})}\frac{E(ax^{1-i-j+k_i+k_j}
q^{\la_i+\la_j-k_i-k_j+1})}
{E(ax^{2-i-j}q)}\right.\hfill\\
\times\left.\frac{(aqx^{3-i-j};q)_{\la_i+\la_j-k_i-k_j}
(x^{j-i+1};q)_{\la_i-\la_j-k_i+k_j}}
{(aq^2x^{1-i-j};q)_{\la_i+\la_j-k_i-k_j}
(qx^{j-i-1};q)_{\la_i-\la_j-k_i+k_j}}\right\}\hfill\\
\times\frac{(aqx^{1-n},b,c,d,eq,q^{-N};q,x)_{\la-k}}
{(qx^{n-1},aq^2/b,aq^2/c,aq^2/d,
aq/e,aq^{N+2};q,x)_{\la-k}}.\hfill\end{multline*}

We will identify the  sum with respect to $k$  as a case of \eqref{w2}
with $q$ replaced by $x^{-1}$. 
Since $k_i\in\{0,1\}$,  we can write
$$(ex^{1-i}q^{\la_i-k_i-N};x)_{k_i}=(ex^{1-i}q^{\la_i-1-N};x^{-1})_{k_i},$$
\begin{equation*}\begin{split}&\quad\frac{(aqx^{1-n},eq,q^{-N};q,x)_{\la-k}}
{(qx^{n-1},aq/e,aq^{N+2};q,x)_{\la-k}}\\
&=\frac{(aqx^{1-n},eq,q^{-N};q,x)_{\la}}
{(qx^{n-1},aq/e,aq^{N+2};q,x)_{\la}}
\prod_{i=1}^n\frac{(x^{n-i}q^{\la_i},
ax^{1-i}q^{\la_i}/e,ax^{1-i}q^{\la_i+N+1};x^{-1})_{k_i}}
{(ax^{2-n-i}q^{\la_i},ex^{1-i}q^{\la_i},
x^{1-i}q^{\la_i-N-1};x^{-1})_{k_i}},\end{split}\end{equation*}
$$(b;q,x)_{\la-k}\prod_{i=1}^n(bx^{1-i}q^{\la_i-k_i};x)_{k_i}=(b;q,x)_\la,$$
$$\frac 1{(aq/b;x^{-1})_n(aq^2/b;q,x)_{\la-k}}\prod_{i=1}^n
\frac{E(ax^{1-i}q^{\la_i-k_i+1}/b)}{(ax^{1-i}q^{\la_i-k_i+1}/b;x)_{k_i}}=
\frac 1{(aq/b;q,x)_\la}$$
and similarly with $b$ replaced by $c$ and $d$.
Using the reflection formula $E(x)=-xE(1/x)$ and recalling \eqref{bc}, we have
$$\frac{E(ax^{i-n}q^{k_i-\la_i+1}/bcd)}
{(ax^{1-i}q^{\la_i-k_i+N+1}/e;x)_{k_i}}=q^{k_i}
\frac{E(ex^{i-1}q^{-\la_i-N-1}/a)}{(ax^{1-i}q^{\la_i+N+1}/e;x^{-1})_{k_i}}.$$
Considering the four cases $k_i$, $k_j=0,1$ separately, we
 find that the factor in curly brackets may be written as
$$\frac{E(x^{j-i+k_j-k_i}q^{\la_i-\la_j},ax^{3-i-j-k_i-k_j}q^{\la_i+\la_j},
aqx^{1-i-j})}{E(x^{j-i},ax^{3-i-j},aqx^{2-i-j})}
\frac{(ax^{3-i-j};q)_{\la_i+\la_j}(x^{j-i+1};q)_{\la_i-\la_j}}
{(aqx^{1-i-j};q)_{\la_i+\la_j}(qx^{j-i-1};q)_{\la_i-\la_j}}.
$$
Finally, by \eqref{dp}, we have
$$(aq;x^{-1})_n\prod_{i=1}^n \frac 1{E(aqx^{2(1-i)})}\prod_{1\leq i<j\leq n}
\frac{E(aqx^{1-i-j})}{E(aqx^{2-i-j})}=1.$$

These simplifications lead to
\begin{multline*}R=\frac{1}
{(aq/bcd;x^{-1})_n}\sum_{\la}\prod_{i=1}^n \left(E(ex^{i-1}q^{-\la_i-N-1}/a)
q^{\la_i}x^{2(i-1)\la_i}\right)\\
\begin{split}&\times\prod_{1\leq i<j\leq n}
\left(\frac{E(x^{j-i}q^{\la_i-\la_j},ax^{2-i-j}q^{\la_i+\la_j})}
{E(x^{j-i},ax^{3-i-j})}
\frac{(ax^{3-i-j};q)_{\la_i+\la_j}(x^{j-i+1};q)_{\la_i-\la_j}}
{(aqx^{1-i-j};q)_{\la_i+\la_j}(qx^{j-i-1};q)_{\la_i-\la_j}}\right)\\
&\times\frac{(b,c,d,aqx^{1-n},eq,q^{-N};q,x)_\la}
{(qx^{n-1},aq/b,aq/c,aq/d,aq/e,aq^{N+2};q,x)_\la}\\
&\times\sum_{k} \prod_{1\leq i<j\leq n}
\frac{E(x^{j-i+k_j-k_i}q^{\la_i-\la_j})}
{E(x^{j-i}q^{\la_i-\la_j})}
\frac{E(ax^{3-i-j-k_i-k_j}q^{\la_i+\la_j})}{E(ax^{2-i-j}q^{\la_i+\la_j})}\\
&\times\prod_{i=1}^n\frac{(x^{n-i}q^{\la_i},
ax^{1-i}q^{\la_i}/e,ax^{1-i}q^{\la_i+N+1},ex^{1-i}q^{\la_i-N-1}
;x^{-1})_{k_i}(-1)^{k_i}
x^{-(i-1)k_i}}
{(ax^{2-n-i}q^{\la_i},ex^{1-i}q^{\la_i},
x^{1-i}q^{\la_i-N-1},ax^{1-i}q^{\la_i+N+1}/e;x^{-1})_{k_i}}.
\end{split}\end{multline*}
The sum in $k$ is the left-hand side of \eqref{w2} with
$$(a,b,c,d,e,x_i,q)\mapsto(ax, x^{n-1}, a/e, aq^{N+1}, eq^{-N-1}, 
x^{1-i}q^{\la_i}, x^{-1}),$$
 and thus equals
\begin{equation*}\begin{split}&\quad
\prod_{i=1}^n\frac{E(ax^{2(1-i)}q^{2\la_i},ex^{i-n},x^{i-n}q^{-N-1},
ex^{i-1}q^{-N-1}/a)}
{E(ax^{2-i-n}q^{\la_i},
ex^{1-i}q^{\la_i},x^{1-i}q^{\la_i-N-1},ex^{i-1}q^{-\la_i-N-1}/a)}\\
&=\frac{(aq/bcd;x^{-1})_n(e,q^{-N-1};q,x)_\la}{(eq,q^{-N};q,x)_\la}
\prod_{i=1}^n\frac{E(ax^{2(1-i)}q^{2\la_i})}
{E(ax^{2-i-n}q^{\la_i},ex^{i-1}q^{-\la_i-N-1}/a)}.
\end{split}\end{equation*}
Finally, we use \eqref{dp} to write
\begin{equation*}\begin{split}&\quad(aqx^{1-n};q,x)_\la\prod_{i=1}^n
\frac 1{E(ax^{2-i-n}q^{\la_i})}
\prod_{1\leq i<j\leq n} \frac 1{E(ax^{3-i-j})}\\
&=(ax^{1-n};q,x)_\la\prod_{i=1}^n\frac 1{E(ax^{2(1-i)})}
\prod_{1\leq i<j\leq n} \frac 1{E(ax^{2-i-j})}.
\end{split}\end{equation*}
Putting all this together we find that $R$ equals
\begin{multline*}\sum_{\la}\prod_{i=1}^n 
\left(\frac{E(ax^{2(1-i)}q^{2\la_i})}{E(ax^{2(1-i)})}\,
q^{\la_i}x^{2(i-1)\la_i}\right)
\frac{(b,c,d,ax^{1-n},e,q^{-N-1};q,x)_\la}
{(qx^{n-1},aq/b,aq/c,aq/d,aq/e,aq^{N+2};q,x)_\la}\hfill\\
\times\prod_{1\leq i<j\leq n}\left(\frac{E(x^{j-i}q^{\la_i-\la_j})}
{E(x^{j-i})}\frac{E(ax^{2-i-j}q^{\la_i+\la_j})}{E(ax^{2-i-j})}
\frac{(ax^{3-i-j};q)_{\la_i+\la_j}(x^{j-i+1};q)_{\la_i-\la_j}}
{(aqx^{1-i-j};q)_{\la_i+\la_j}(qx^{j-i-1};q)_{\la_i-\la_j}}\right),
\hfill
\end{multline*}
which is indeed the left-hand side of \eqref{w} with $N$ replaced by $N+1$.
This completes the proof of Theorem \ref{t}. 

\section{Duality}

In this section we prove a duality property for sums
of the type occurring in Theorem \ref{t}. To state the result, we
use the notation
\begin{multline*}
{}_{r+1}\Om_r^{(n)}(a;b_1,\dots,b_{r-3},q^{-N};q,x)\\
=\sum_{\la\in\La_{nN}}\prod_{i=1}^n 
\left(\frac{E(ax^{2(1-i)}q^{2\la_i})}{E(ax^{2(1-i)})}
\,q^{\la_i}x^{2(i-1)\la_i}\right)\prod_{1\leq i<j\leq n}\left(\frac{E(x^{j-i}
q^{\la_i-\la_j})}{E(x^{j-i})}\frac{E(ax^{2-i-j}q^{\la_i+\la_j})}
{E(ax^{2-i-j})}\right.\\
\times\left.\frac{(ax^{3-i-j};q)_{\la_i+\la_j}
(x^{j-i+1};q)_{\la_i-\la_j}}{(aqx^{1-i-j};q)_{\la_i+\la_j}
(qx^{j-i-1};q)_{\la_i-\la_j}}\right)
\frac{(ax^{1-n},b_1,\dots,b_{r-3},q^{-N};q,x)_\la}{(qx^{n-1},aq/b_1,\dots,
aq/b_{r-3},aq^{N+1};q,x)_\la}.
\end{multline*}
It is natural to assume 
 the balanced condition
$(aq)^{r-3}=(x^{n-1}q^{1-N}\prod_{i}b_i)^2$, though we do not need it
to prove the following proposition. 

\begin{proposition}\label{p}
One has
$${}_{r+1}\Om_r^{(n)}(a;b_1,\dots,b_{r-3},q^{-N};q,x)
={}_{r+1}\Om_r^{(N)}(aqx;b_1,\dots,b_{r-3},x^{n};x^{-1},q^{-1}).$$
\end{proposition}

In fact, this  holds as a termwise symmetry between the two sums, 
the change of summation variable $\La_{nN}\rightarrow\La_{Nn}$
 being conjugation of
partitions. Let us write $\la'$ for the conjugate of a partition $\la$.
Note that, since we consider partitions into
non-negative parts, $\la'$ depends not only on the Young diagram
of $\la$ but also on the choice of $n$ and $N$. 
For instance, $(3,2,0)\in\La_{33}$ and $(3,2,0)\in\La_{34}$
has conjugate $(2,2,1)$ and $(2,2,1,0)$, respectively.

To prove Proposition \ref{p}, we  observe that since clearly
 $(b;q,x)_\la=(b;x^{-1},q^{-1})_{\la'}$, it is enough to show 
that, for $\la\in\La_{nN}$, the two quantities
$$A_\la=
\prod_{1\leq i\leq j\leq n}\frac{E(ax^{2-i-j}q^{\la_i+\la_j})}{E(ax^{2-i-j})}
\prod_{1\leq i<j\leq n}
\frac{(ax^{3-i-j};q)_{\la_i+\la_j}}{(aqx^{1-i-j};q)_{\la_i+\la_j}}
\frac{(ax^{1-n};q,x)_\la}{(aq^{N+1};q,x)_\la}$$
and 
\begin{equation*}B_\la=
\prod_{i=1}^n q^{\la_i}x^{2(i-1)\la_i}\prod_{1\leq i<j\leq n}
\frac{E(x^{j-i}q^{\la_i-\la_j})}{E(x^{j-i})}
\frac{(x^{j-i+1};q)_{\la_i-\la_j}}{(qx^{j-i-1};q)_{\la_i-\la_j}}
\frac{(q^{-N};q,x)_\la}{(qx^{n-1};q,x)_\la}
\end{equation*}
are invariant under the transformation 
$(a,q,x,n,N,\la)\mapsto(aqx,x^{-1},q^{-1},N,n,\la')$.

We prove the invariance of $A_\la$, the  case of $B_\la$ being 
similar. We fix $n$ and $N$ and proceed by induction on the 
number of boxes in the Young diagram of $\la$, 
starting from the trivial case of zero boxes.
Suppose that  the invariance holds for a fixed partition 
$\la$.
We will show that it also holds for any partition 
$\la^+$ obtained by adding a box to the Young diagram of $\la$. 
There exist $k$ and $l$ with $1\leq k\leq n$, $1\leq l\leq N$ such that
$\la_j^+=\la_j$ for $j\neq k$, $\la_k=l-1$ and $\la_k^+=l$.
After straight-forward simplifications, we may write
$$\frac{A_{\la^+}}{A_\la}
=\frac{E(ax^{2-2k}q^{2l},ax^{1-2k}q^{2l-1},ax^{2-n-k}q^{l-1})}
{E(ax^{2-2k}q^{2l-1},ax^{3-2k}q^{2l-2},ax^{1-k}q^{l+N})}
\prod_{i=1}^n\frac{E(ax^{2-i-k}q^{\la_i+l},ax^{3-i-k}q^{\la_i+l-1})}
{E(ax^{2-i-k}q^{\la_i+l-1},ax^{1-i-k}q^{\la_i+l})}.$$
Next we observe that
\begin{equation*}\begin{split}
&\quad\prod_{i=1}^n\frac{E(bx^{1-i}q^{\la_i})}{E(bx^{-i}q^{\la_i})}
=\prod_{i=1}^{\la_N'}\frac{E(bx^{1-i}q^N)}{E(bx^{-i}q^N)}
\prod_{i=\la_N'+1}^{\la_{N-1}'}\frac{E(bx^{1-i}q^{N-1})}{E(bx^{-i}q^{N-1})}
\times\dotsm\times\prod_{i=\la_1'+1}^{n}\frac{E(bx^{1-i})}{E(bx^{-i})}\\
&=\frac{E(bq^N)}{E(bx^{-\la_N'}q^N)}
\frac{E(bx^{-\la_N'}q^{N-1})}{E(bx^{-\la_{N-1}'}q^{N-1})}\dotsm
\frac{ E(bx^{-\la_1'})}{E(bx^{-n})}
=\frac{E(bq^N)}{E(bx^{-n})}\prod_{i=1}^N
\frac{E(bx^{-\la_i'}q^{i-1})}{E(bx^{-\la_i'}q^i)},
\end{split}\end{equation*}
which gives
$$\frac{A_{\la^+}}{A_\la}
=\frac{E(ax^{2-2k}q^{2l},ax^{1-2k}q^{2l-1},ax^{2-k}q^{l+N-1})}
{E(ax^{2-2k}q^{2l-1},ax^{3-2k}q^{2l-2},ax^{1-k-n}q^{l})}
\prod_{i=1}^N\frac{E(ax^{1-k-\la_i'}q^{l+i-1},ax^{2-k-\la_i'}q^{l+i-2})}
{E(ax^{1-k-\la_i'}q^{l+i},ax^{2-k-\la_i'}q^{l+i-1})}.$$
This agrees with the expression  obtained from the previous one by 
 substituting $(a,q,x,n,N,\la,k,l)\mapsto(aqx,x^{-1},q^{-1},N,n,\la',l,k)$.
Thus the  invariance of $A_\la$ implies that of $A_{\la^+}$.

\end{document}